\documentclass[11pt]{article}
\usepackage{latexsym,amssymb,amsmath,amscd,amsfonts,mathrsfs,texdraw, colordvi}

\newtheorem{theorem}{Theorem}[section]

\numberwithin{equation}{section} \numberwithin{figure}{section}
\newtheorem{proposition}[theorem]{Proposition}
\newtheorem{corollary}[theorem]{Corollary}
\def\qed{\hfill \rule{4pt}{7pt}}

\def\pf{\noindent {\it Proof.} }

\begin{document}
%-----------------------------------------------------------------------

\begin{center}
{\large {\bf Oscillating Rim Hook Tableaux and Colored Matchings}}

\vskip 6mm

{\small William Y.C. Chen$^1$
and Peter L. Guo$^2$\\[%
2mm] Center for Combinatorics, LPMC-TJKLC\\
Nankai University, Tianjin 300071,
P.R. China \\[3mm]
$^1$chen@nankai.edu.cn,
$^2$lguo@cfc.nankai.edu.cn \\[0pt%
] }
\end{center}

%------------------------------------------------------------------------
\begin{abstract}
Motivated by the question of finding a type $B$ analogue of the
bijection between oscillating tableaux
 and matchings, we find a correspondence between
 oscillating $m$-rim hook tableaux and $m$-colored
 matchings, where $m$ is a positive integer. An oscillating $m$-rim hook tableau is
 defined as a sequence $(\lambda^0,\lambda^1,\ldots,\lambda^{2n})$
  of Young diagrams
 starting with the empty shape and ending with the
 empty shape such that $\lambda^{i}$
 is obtained from $\lambda^{i-1}$ by adding
 an $m$-rim hook or by deleting an $m$-rim hook.
Our bijection relies on the generalized Schensted algorithm due to
White. An oscillating $2$-rim hook tableau is also called an
oscillating domino tableau. When we restrict our attention to two
column oscillating domino tableaux of length $2n$, we are led to a
bijection between such tableaux and noncrossing $2$-colored
matchings on $\{1, 2,\ldots, 2n\}$, which are counted by the product
$C_nC_{n+1}$ of two consecutive Catalan numbers. A 2-colored
matching
 is  noncrossing
 if there are no two arcs of the same color that
 are intersecting.
We   show that  oscillating domino tableaux with at most two columns are in
one-to-one correspondence with Dyck path packings. A Dyck path
packing of length $2n$ is a pair $(D,E)$, where
 $D$ is a Dyck path of length $2n$,
 and $E$ is a dispersed Dyck path of length $2n$
  that is  weakly covered by
$D$. So we deduce  that
 Dyck path packings of length $2n$ are counted
by  $C_nC_{n+1}$.

 \end{abstract}

\vskip 3mm


%-------------------------------------------------------------
\allowdisplaybreaks

\section{Introduction}

The objective of this paper is to provide a rim hook analogue of the
correspondence between
 oscillating tableaux and matchings given by
Chen, Deng, Du, Stanley and Yan \cite{Chen}.
 We show that there is
a one-to-one correspondence between oscillating $m$-rim hook
tableaux and $m$-colored matchings. The construction of our
bijection relies on the generalized Schensted
 algorithm for rim
hook tableaux introduced by White \cite{White}.

\begin{figure}[h,t]
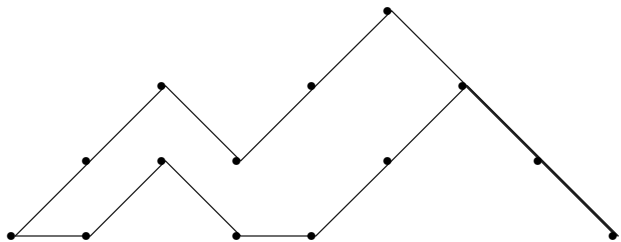

\centertexdraw{ \drawdim mm \linewd 0.2 \setgray 0.1
%\move(0
%0)\arrowheadtype t:V \avec(100 0)\move(0 0)\arrowheadtype t:V
%\avec(0 40)
\move(0 0)\lvec(10 10)\lvec(20 20)\lvec(30 10)\lvec(40 20)\lvec(50
30)\lvec(60 20)\lvec(70 10)\lvec(80 0)
%\move(0
%0)\arrowheadtype t:V \avec(0 40)
\move(0 0)\lvec(10 0)\lvec(20 10)\lvec(30 0)\lvec(40 0)\lvec(50
10)\lvec(60 20)

\textref h:C v:C \htext(0 0){\circle*{2.5}}\textref h:C v:C
\htext(10 0){\circle*{2.5}}\textref h:C v:C \htext(10
10){\circle*{2.5}}\textref h:C v:C \htext(20
10){\circle*{2.5}}\textref h:C v:C \htext(20
20){\circle*{2.5}}\textref h:C v:C \htext(30
0){\circle*{2.5}}\textref h:C v:C \htext(30
10){\circle*{2.5}}\textref h:C v:C \htext(40
0){\circle*{2.5}}\textref h:C v:C \htext(40 20){\circle*{2.5}}
\textref h:C v:C \htext(50 10){\circle*{2.5}}\textref h:C v:C
\htext(50 30){\circle*{2.5}}\textref h:C v:C \htext(60
20){\circle*{2.5}}\textref h:C v:C \htext(70 10){\circle*{2.5}}
\textref h:C v:C \htext(80 0){\circle*{2.5}}

\move(60.15 20)\lvec(80.15 0)\move(60 19.85)\lvec(79.85 0)

} \caption{A Dyck path packing.}\label{pair}
\end{figure}

We shall pay special attention to the case of oscillating domino
(2-rim hook) tableaux with at most two columns. In this case, our
main result reduces to a bijection between oscillating domino
tableaux with at most two columns and noncrossing $2$-colored
matchings. Bear in mind that, a noncrossing $2$-colored
 matching is not meant to be a noncrossing matching with
two colors, but a matching that does not contain crossing edges that
are of the same color. On the other hand, we find a correspondence
between oscillating domino tableaux with at most two columns and
Dyck path packings.
 A Dyck path packing can be viewed as a
 dispersed Dyck path $E$ weakly covered
 by a Dyck path $D$, where a dispersed Dyck path
 is defined as non-overlapping Dyck paths connected by
 some horizontal steps on the $x$-axis. See
 Figure \ref{pair} for an illustration.

So we are led to a bijection between Dyck path packings and
noncrossing $2$-colored matchings. It is easy to check that
noncrossing $2$-colored matchings with $2n$ vertices are counted by
the product of two Catalan numbers, that is, $C_nC_{n+1}$, where
\[ C_n=\frac{1}{n+1}{2n \choose n}.\]

Notice that there are several combinatorial  objects that are
enumerated by the number  $C_nC_{n+1 }$, such as  walks within the
first quadrant of $\mathbb{Z}^3$ starting at $(0,0,0)$
 and consisting of $n$ steps taken from $\{(-1,0,0), (0,-1,1), (0,1,0), (1,0,
-1)\}$, see  Bostan and Kauers \cite{Bostan},  walks within the
first quadrant of $\mathbb{Z}^2$ starting at $(0,0)$, ending on the
$x$-axis and consisting of $2n$ steps taken from $\{(-1,0), (-1,1),
(1,-1), (1,0)\}$, see Bousquet-M\'elou and Mishna \cite{Melou},
alternating
 Baxter permutations
of length $2n+1$, see Cori, Dulucq and Viennot \cite{Cori}, and
 walks within the
first quadrant of $\mathbb{Z}^2$ starting and ending at $(0,0)$ and
consisting of $2n$ steps taken from $\{(-1,0), (0,-1), (0,1),
(1,0)\}$, see Guy \cite{Guy}.  We also find a correspondence between
noncrossing $2$-colored matchings and Guy's walks. It would be
interesting to establish further connections between packings of
Dyck paths and other combinatorial structures.

\begin{figure}[h,t]
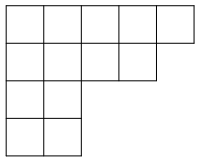

\centertexdraw{ \drawdim mm \linewd 0.1 \setgray 0 \move(0
0)\lvec(10 0)\move(0 5)\lvec(10 5)\move(0 10)\lvec(20 10) \move(0
15)\lvec(25 15)\move(0 20)\lvec(25 20) \move(0 0)\lvec(0 20)\move(5
0)\lvec(5 20)\move(10 0)\lvec(10 20)\move(15 10)\lvec(15 20)\move(20
10)\lvec(20 20)\move(25 15)\lvec(25 20)

} \caption{The Young diagram of $(5,4,2,2)$}\label{Ferrers1}
\end{figure}

Let us give an overview of some definitions.
 A partition  of an integer $n$  is
a sequence $(\lambda_1,\lambda_2,\ldots,\lambda_\ell)$ of
nonincreasing positive integers such that
$\lambda_1+\lambda_2+\cdots+\lambda_\ell=n$. We can also represent a
partition  by its Young diagram, i.e., a left-justified array of
cells  (or, squares) with $\lambda_i$ cells in  row $i$ for $1\leq i
\leq \ell$. For example,  Figure \ref{Ferrers1} is the Young diagram
of the partition $(5,4,2,2)$.

\begin{figure}[h,t]
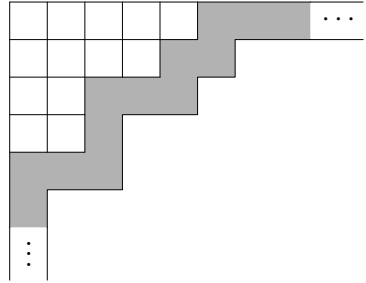

\centertexdraw{ \drawdim mm \linewd 0.1 \setgray 0
 \move(40
15)\lvec(30 15)\lvec(30 10)\lvec(25 10)\lvec(25 5)\lvec(15
5)\lvec(15 -5)\lvec(5 -5)\lvec(5 -10) \lvec(0 -10)\lvec(0 0)\lvec(10
0)\lvec(10 10)\lvec(20 10)\lvec(20 15)\lvec(25 15)\lvec(25
20)\lvec(40 20) \ifill f:0.7

\move(0 0)\lvec(10 0)\move(0 5)\lvec(10 5)\move(0 10)\lvec(20 10)
\move(0 15)\lvec(25 15)\move(0 20)\lvec(40 20) \move(0 -10)\lvec(0
20)\move(5 0)\lvec(5 20)\move(10 0)\lvec(10 20)\move(15 10)\lvec(15
20)\move(20 10)\lvec(20 20)\move(25 15)\lvec(25 20) \move(40
15)\lvec(30 15)\lvec(30 10)\lvec(25 10)\lvec(25 5)\lvec(15
5)\lvec(15 -5)\lvec(5 -5)\lvec(5 -10)

\move(40 20)\lvec(47 20)\move(40 15)\lvec(47 15)

\move(0 -10)\lvec(0 -17)\move(5 -10)\lvec(5 -17)

\textref h:C v:C \htext(44 17.5){$\cdots$}\textref h:C v:C
\htext(2.5 -12.5){$\vdots$}

} \caption{The outside border of $(5,4,2,2)$}\label{outerborder}
\end{figure}

To define  the rim hooks of a partition, we note that the outside
border of a partition $\lambda$ is a collection of cells not in
$\lambda$ but immediately bellow or to the right of $\lambda$; or in
the first row and to the right of $\lambda$; or in the first column
and bellow $\lambda$. For example, in Figure \ref{outerborder} the
shaded area illustrates the outside border of $\lambda=(5,4,2,2)$.

\begin{figure}[h,t]
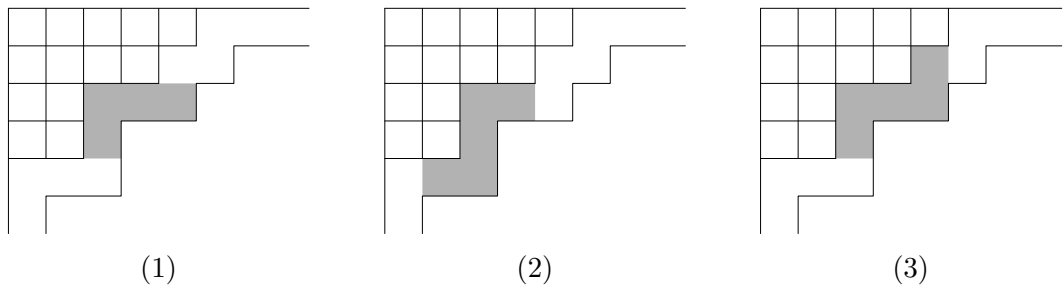

\centertexdraw{ \drawdim mm \linewd 0.1 \setgray 0 \move(25
25)\lvec(25 10)\lvec(10 10)\lvec(10 0)\lvec(15 0)\lvec(15 5)\lvec(25
5)\ifill f:0.7

\move(70 10)\lvec(60 10)\lvec(60 0)\lvec(55 0)\lvec(55 -5)\lvec(65
-5)\lvec(65 5)\lvec(70 5)\ifill f:0.7

\move(125 5)\lvec(125 15)\lvec(120 15)\lvec(120 10)\lvec(110
10)\lvec(110 0)\lvec(115 0)\lvec(115 5)\ifill f:0.7

 \move(0 0)\lvec(10 0)\move(0 5)\lvec(10 5)\move(0 10)\lvec(20
10) \move(0 15)\lvec(25 15)\move(0 20)\lvec(40 20) \move(0
-10)\lvec(0 20)\move(5 0)\lvec(5 20)\move(10 0)\lvec(10 20)\move(15
10)\lvec(15 20)\move(20 10)\lvec(20 20)\move(25 15)\lvec(25 20)
\move(40 15)\lvec(30 15)\lvec(30 10)\lvec(25 10)\lvec(25 5)\lvec(15
5)\lvec(15 -5)\lvec(5 -5)\lvec(5 -10)

\move(50 0)\lvec(60 0)\move(50 5)\lvec(60 5)\move(50 10)\lvec(70 10)
\move(50 15)\lvec(75 15)\move(50 20)\lvec(90 20) \move(50
-10)\lvec(50 20)\move(55 0)\lvec(55 20)\move(60 0)\lvec(60
20)\move(65 10)\lvec(65 20)\move(70 10)\lvec(70 20)\move(75
15)\lvec(75 20) \move(90 15)\lvec(80 15)\lvec(80 10)\lvec(75
10)\lvec(75 5)\lvec(65 5)\lvec(65 -5)\lvec(55 -5)\lvec(55 -10)

\move(100 0)\lvec(110 0)\move(100 5)\lvec(110 5)\move(100
10)\lvec(120 10) \move(100 15)\lvec(125 15)\move(100 20)\lvec(140
20) \move(100 -10)\lvec(100 20)\move(105 0)\lvec(105 20)\move(110
0)\lvec(110 20)\move(115 10)\lvec(115 20)\move(120 10)\lvec(120
20)\move(125 15)\lvec(125 20) \move(140 15)\lvec(130 15)\lvec(130
10)\lvec(125 10)\lvec(125 5)\lvec(115 5)\lvec(115 -5)\lvec(105
-5)\lvec(105 -10)\textref h:C v:C \htext(20 -15){$(1)$} \textref h:C
v:C \htext(70 -15){$(2)$}\textref h:C v:C \htext(120 -15){$(3)$} }
\caption{Contiguous cells in the outside border}\label{illegal}
\end{figure}

Let $\alpha$ be a set of contiguous cells in the outside border of
$\lambda$. We say that $\alpha$ is a rim hook outside $\lambda$
 if the
shape $\mu=\lambda\cup\alpha$ is the Young diagram of a partition.
For example, in Figure \ref{illegal}, it can be seen that among the
three sets of contiguous cells in
 the outside border of
$\lambda$,  there is only one rim hook outside $\lambda$, which is
the third diagram.
If $\alpha$ has $m$ cells, then
 we write $|\alpha|=m$ and call $\alpha$ an $m$-rim hook outside $\lambda$.
If $\alpha$ is a rim hook outside $\lambda$ and
$\mu=\lambda\cup\alpha$, then we call $\alpha$ an outer rim hook of
$\mu$, and we write $\mu-\alpha$ to mean $\lambda$.

\begin{figure}[h,t]
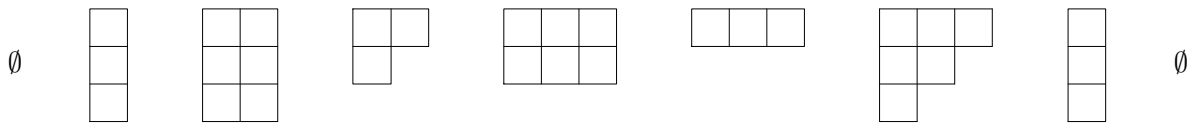

\centertexdraw{ \drawdim mm \linewd 0.1 \setgray 0

\move(10 0)\lvec(15 0)\lvec(15 15)\lvec(10 15)\lvec(10 0)\move(10
5)\lvec(15 5)\move(10 10)\lvec(15 10)

\move(25 0)\lvec(35 0)\lvec(35 15)\lvec(25 15)\lvec(25 0)\move(25
5)\lvec(35 5)\move(30 0)\lvec(30 15)\move(25 10)\lvec(35 10)

\move(45 15)\lvec(55 15)\lvec(55 10)\lvec(50 10)\lvec(50 5)\lvec(45
5)\lvec(45 15)\move(45 10)\lvec(50 10)\lvec(50 15)

\move(65 5)\lvec(80 5)\lvec(80 15)\lvec(65 15)\lvec(65 5) \move(65
10)\lvec(80 10)\move(70 5)\lvec(70 15)\move(75 5)\lvec(75 15)

\move(90 10)\lvec(105 10)\lvec(105 15)\lvec(90 15)\lvec(90 10)
\move(95 10)\lvec(95 15)\move(100 10)\lvec(100 15)

\move(115 15)\lvec(130 15)\move(115 10)\lvec(130 10)\move(115
5)\lvec(125 5)\move(115 0)\lvec(120 0)\move(115 15)\lvec(115 0)
\move(120 15)\lvec(120 0)\move(125 15)\lvec(125 5)\move(130
15)\lvec(130 10)

\move(140 0)\lvec(145 0)\lvec(145 15)\lvec(140 15)\lvec(140 0)
\move(140 5)\lvec(145 5)\move(140 10)\lvec(145 10)

 \textref h:C v:C \htext(155
8){$\emptyset$}\textref h:C v:C \htext(0 8){$\emptyset$}

} \caption{An oscillating 3-rim hook  tableau}\label{oscillating1}
\end{figure}

We can now introduce the notion of
  oscillating $m$-rim hook tableaux.
An oscillating $m$-rim hook tableau
 of length $2n$
can be defined as a sequence
$\lambda=(\lambda^0,\lambda^1,\ldots,\lambda^{2n})$ of Young
diagrams such that $\lambda^0=\lambda^{2n}=\emptyset$,
 and for
$1\leq i\leq 2n$, $\lambda^i$ is obtained
 from $\lambda^{i-1}$
either by adding an $m$-rim hook
 outside $\lambda^{i-1}$ or by
deleting  an  outer $m$-rim hook  of $\lambda^{i-1}$. For example,
Figure \ref{oscillating1} is an illustration of an oscillating 3-rim
hook tableau.

When $m=1$, an oscillating $m$-rim hook tableau
 is an ordinary oscillating tableau.
 An oscillating $2$-rim hook tableau
will be also called an oscillating domino tableau. We shall use
$r(\lambda)$ (resp., $c(\lambda)$) to denote the maximum number of
rows (resp., columns) of the shape $\lambda^i$ appearing in
$\lambda$ for $0\leq i \leq 2n$.

The main objective of this paper is to show that there is a
one-to-one correspondence between oscillating $m$-rim hook tableaux
and $m$-colored matchings. An $m$-colored matching $M$ on $[2n]$ is
a matching on $[2n]$ with each arc assigned one of  $m$  colors,
say, ${c}_1,{c}_2,\ldots,{c}_m$. For example,
 Figure 1.6 gives a 2-colored matching,
  where we use solid lines to
represent  arcs assigned the color ${c}_1$,
 and dotted lines to
represent arcs assigned the color ${c}_2$.
\begin{figure}[h,t]
\setlength{\unitlength}{0.4mm}
\begin{center}
\begin{picture}(260,50)
\put(0,10){\circle*{2}}\put(20,10){\circle*{2}}\put(40,10){\circle*{2}}
\put(60,10){\circle*{2}}\put(80,10){\circle*{2}}\put(100,10){\circle*{2}}
\put(120,10){\circle*{2}}\put(140,10){\circle*{2}}
\put(160,10){\circle*{2}}\put(180,10){\circle*{2}}
\put(200,10){\circle*{2}}\put(220,10){\circle*{2}}
\put(240,10){\circle*{2}}\put(260,10){\circle*{2}}
\qbezier(0,10)(60,60)(120,10)\qbezier[60](20,10)(60,50)(100,10)
\qbezier(40,10)(130,75)(220,10)\qbezier(60,10)(130,65)(200,10)
\qbezier[120](80,10)(180,80)(260,10)\qbezier[70](140,10)(190,50)(240,10)
\qbezier(160,10)(170,20)(180,10)

\put(-2,0){1}\put(17.5,0){2}\put(38,0){3}\put(58,0){4}
\put(78,0){5}\put(98,0){6}\put(118,0){7}
\put(-2,0){1}\put(138.5,0){8}\put(158.5,0){9} \put(175.5,0){10}
\put(195.5,0){11}\put(215.5,0){12}\put(235.5,0){13}\put(255.5,0){14}
\end{picture}\caption{A 2-colored matching}
\end{center}\label{mn}
\end{figure}
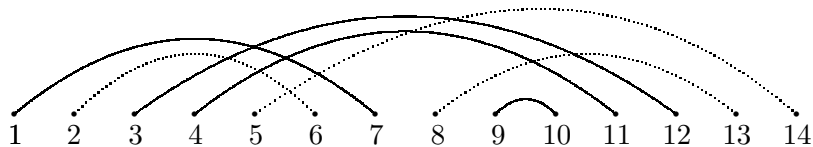A $k$-crossing of $M$ is a $k$-subset
$\{(i_1,j_1),(i_2,j_2),\ldots,(i_k,j_k)\}$ of arcs of the same color
such that $i_1<i_2<\cdots <i_k<j_1<j_2<\cdots<j_k$. Similarly, we
define a $k$-nesting of $M$ as a $k$-subset
$\{(i_1,j_1),(i_2,j_2),\ldots,(i_k,j_k)\}$ of arcs of the same color
such that $i_1<i_2<\cdots <i_k<j_k<\cdots<j_2<j_1$. Denote by
cr($M$) (resp., ne($M$)) the maximal number
 $k$ such that $M$ has a
$k$-crossing (resp., $k$-nesting). We say that $M$ is
$k$-noncrossing (resp., $k$-nonnesting) if $M$ has no $k$-crossing
(resp., $k$-nesting). A 2-noncrossing (or, 2-nonnesting) $m$-colored
matching is  called a noncrossing (or, nonnesting) $m$-colored
matching.

Our bijection can be used to characterize the crossing number and
the nesting number of an $m$-colored matching in terms of the
maximum number of columns and the maximum number of rows of shapes
 in the
corresponding oscillating $m$-rim hook tableau.
 The construction of
our bijection is based on the generalized Schensted algorithm for
rim hook tableaux due to White \cite{White}.

This paper is organized as follows.
 We shall give a
brief review  of White's algorithm in  Section 2. Based on this
algorithm, we give a description of the bijection between
oscillating $m$-rim hook tableaux and $m$-colored matchings in
Section 3. Section 4 is concerned with oscillating domino tableaux
with at most two columns and noncrossing $2$-colored matchings.
 We show that
such tableaux are in one-to-one correspondence with
 Dyck path packings. We also give
a bijection between noncrossing $2$-colored matchings and Guy's
walks.

%-------------------------------------------------------------------------

\section{The generalized Schensted algorithm}

To give a combinatorial proof of the orthogonality of the characters
of the symmetric group $S_n$, White \cite{White} extended the
ordinary  Schensted algorithm \cite{Schensted} to rim hook tableaux.
For our purpose, we shall be
 concerned with only  a
special  case of White's construction
when all rim hooks and
hooks are restricted to $m$-rim hooks and $m$-hooks.  This version of
White's algorithm has been further studied by
 Stanton and White
\cite{Stanton}. As remarked  by White \cite{White}, when $m=2$,
White's algorithm reduces to a one-to-one correspondence between
elements of the hyperoctahedral group and pairs of domino tableaux
of the same shape, which was first obtained  by Lusztig. We shall
adopt the notation and terminology in \cite{Stanton}.

\begin{figure}[h,t]
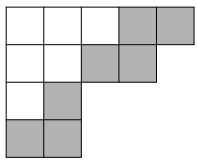

\centertexdraw{ \drawdim mm \linewd 0.1 \setgray 0 \move(0
0)\lvec(10 0)\lvec(10 10)\lvec(20 10)\lvec(20 15)\lvec(25
15)\lvec(25 20)\lvec(15 20)\lvec(15 15)\lvec(10 15)\lvec(10
10)\lvec(5 10)\lvec(5 5)\lvec(0 5) \ifill f:0.7

\move(0 0)\lvec(10 0)\move(0 5)\lvec(10 5)\move(0 10)\lvec(20 10)
\move(0 15)\lvec(25 15)\move(0 20)\lvec(25 20) \move(0 0)\lvec(0
20)\move(5 0)\lvec(5 20)\move(10 0)\lvec(10 20)\move(15 10)\lvec(15
20)\move(20 10)\lvec(20 20)\move(25 15)\lvec(25 20)  } \caption{The
skew diagram of shape $(5,4,2,2)/(3,2,1)$}\label{skewdiagram}
\end{figure}

For two partitions
$\lambda=(\lambda_1,\lambda_2,\ldots)$ and
$\mu=(\mu_1,\mu_2,\ldots)$,
we write $\mu\subseteq \lambda$ if
$\mu_i\leq \lambda_i$ for all $i$.
 If $\mu\subseteq\lambda$, then
the skew diagram of shape $\lambda/\mu$ is
defined as the set of cells obtained from
$\lambda$ by deleting the cells in $\mu$.
For example, the shaded area in Figure
\ref{skewdiagram}  represents a skew diagram.

Let $i_1,i_2,\ldots,i_n$ be $n$ positive integers with
$i_1<i_2<\cdots<i_n$. An $m$-rim hook tableau $P$ of shape $\lambda$
on $\{i_1,i_2,\ldots,i_n\}$ is an assignment of $i_1,i_2,\ldots,i_n$
to the squares of $\lambda$ with  each integer  appearing  exactly
$m$ times such that the set $\alpha$ of squares occupied by $i_n$ is
an outer $m$-rim hook of $\lambda$, and the tableau obtained from
$P$ by deleting the squares occupied by $i_n$ is an $m$-rim hook
tableau on $\{i_1,i_2,\ldots,i_{n-1}\}$.
 We say that
$i_1,i_2,\ldots,i_n$ are the contents of $P$, and write
$content(P)=\{i_1,i_2,\ldots,i_n\}$. Similarly,
we can define   skew
$m$-rim hook tableaux of shape $\lambda/\mu$. Figure \ref{4tableau}
is an illustration of a $4$-rim hook tableau.
\begin{figure}[h,t]
\centertexdraw{ \drawdim mm \linewd 0.1 \setgray 0

\move(0 0)\lvec(15 0)\move(0 5)\lvec(15 5)\move(0 10)\lvec(15 10)
\move(0 15)\lvec(25 15)\move(0 20)\lvec(30 20)\move(0 25)\lvec(30
25)

\move(0 0)\lvec(0 25)\move(5 0)\lvec(5 25)\move(10 0)\lvec(10 25)
\move(15 0)\lvec(15 25)\move(20 15)\lvec(20 25)\move(25 15)\lvec(25
25)\move(30 20)\lvec(30 25)

\textref h:C v:C \htext(2.5 22.5){1}\textref h:C v:C \htext(7.5
22.5){1}\textref h:C v:C \htext(12.5 22.5){1}\textref h:C v:C
\htext(17.5 22.5){1}

\textref h:C v:C \htext(2.5 2.5){2}\textref h:C v:C \htext(2.5
7.5){2}\textref h:C v:C \htext(2.5 12.5){2}\textref h:C v:C
\htext(2.5 17.5){2}

\textref h:C v:C \htext(7.5 7.5){3}\textref h:C v:C \htext(7.5
12.5){3}\textref h:C v:C \htext(7.5 17.5){3}\textref h:C v:C
\htext(12.5 17.5){3}

\textref h:C v:C \htext(17.5 17.5){4}\textref h:C v:C \htext(22.5
17.5){4}\textref h:C v:C \htext(22.5 22.5){4}\textref h:C v:C
\htext(27.5 22.5){4}

\textref h:C v:C \htext(7.5 2.5){5}\textref h:C v:C \htext(12.5
2.5){5}\textref h:C v:C \htext(12.5 7.5){5}\textref h:C v:C
\htext(12.5 12.5){5}

 }
\caption{A $4$-rim hook tableau}\label{4tableau}
\end{figure}

Recall that an $m$-hook is the
 Young diagram corresponding to
a partition $(t,1,1,\ldots,1)$ of $m$,
where $1\leq t\leq m$.
 An $m$-hook
tableau is an $m$-rim hook tableau
whose shape is an $m$-hook.
Figure \ref{hooktableau} gives  a 4-hook tableau
of shape $(2,1,1)$.
\begin{figure}[h,t]
\centertexdraw{ \drawdim mm \linewd 0.1 \setgray 0

\move(0 0)\lvec(5 0)\lvec(5 10)\lvec(10 10)\lvec(10 15)\lvec(0
15)\lvec(0 0)

\move(0 5)\lvec(5 5) \move(0 10)\lvec(5 10)

\move(5 10)\lvec(5 15)

\textref h:C v:C \htext(2.5 2.5){3}\textref h:C v:C \htext(2.5
7.5){3}\textref h:C v:C \htext(2.5 12.5){3}\textref h:C v:C
\htext(7.5 12.5){3}
 }
\caption{A $4$-hook tableau}\label{hooktableau}
\end{figure}

Let $P$ be an $m$-rim hook tableau,
and $H$ be an $m$-hook tableau.
The generalized Schensted algorithm
 is an algorithm to generate
 an  $m$-rim hook tableau by inserting $H$ to $P$.
  When
$m=1$, it reduces to the usual Schensted algorithm.
 To describe the
generalized Schensted algorithm, we need recall
more definitions.

Let $\alpha$ be a set of contiguous
 cells contained  in the outside
border of $\lambda$. The head (resp., tail)
of $\alpha$ is the
upper rightmost (resp., lower leftmost)
square in $\alpha$. The head (resp., tail) of
$\alpha$ is said to be illegal
with respect to  $\lambda$ if the
cell above (resp., to the left) the head (resp., tail)
 is in the outside border of $\lambda$.
 For example, in Figure \ref{illegal} the contiguous
 cells in (1)  have an illegal head, whereas the
 contiguous cells in (2)
have an illegal tail.
Clearly, if $\alpha$ has neither an illegal
head nor  an illegal
tail with respect to $\lambda$, then $\alpha$
is a rim hook outside $\lambda$.
Even though $\alpha$ is not a rim
hook outside $\lambda$, it may be  a
rim hook outside another partition.
We shall say that $\alpha$
is a rim hook if it is a rim
hook outside a certain partition.

\begin{figure}[h,t]
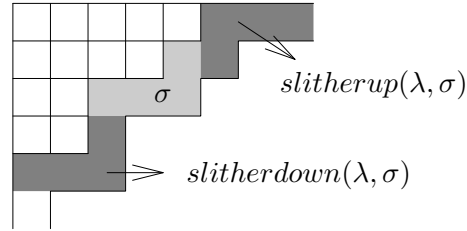

\centertexdraw{ \drawdim mm \linewd 0.1 \setgray 0

\move(10 5)\lvec(25 5)\lvec(25 15)\lvec(20 15)\lvec(20 10)\lvec(10
10)\ifill f:0.8

\move(25 10)\lvec(30 10)\lvec(30 15)\lvec(40 15)\lvec(40 20)\lvec(25
20)\ifill f:0.5

\move(0 -5)\lvec(15 -5)\lvec(15 5)\lvec(10 5)\lvec(10 0)\lvec(0
0)\ifill f:0.5 \textref h:C v:C \htext(20  7.5){$\sigma$}

 \move(0
0)\lvec(10 0)\move(0 5)\lvec(10 5)\move(0 10)\lvec(20 10) \move(0
15)\lvec(25 15)\move(0 20)\lvec(40 20) \move(0 -10)\lvec(0
20)\move(5 0)\lvec(5 20)\move(10 0)\lvec(10 20)\move(15 10)\lvec(15
20)\move(20 10)\lvec(20 20)\move(25 15)\lvec(25 20) \move(40
15)\lvec(30 15)\lvec(30 10)\lvec(25 10)\lvec(25 5)\lvec(15
5)\lvec(15 -5)\lvec(5 -5)\lvec(5 -10)

\move(30 17.5)\arrowheadtype t:V \avec(37.5 12.5)\textref h:C v:C
\htext(48 9){$slitherup(\lambda,\sigma)$} \move(12.5
-2.5)\arrowheadtype t:V \avec(20 -2.5)\textref h:C v:C \htext(38
-3){$slitherdown(\lambda,\sigma)$} }
\caption{$slitherup(\lambda,\sigma)$ and
$slitherdown(\lambda,\sigma)$}\label{slipe}
\end{figure}

Let $\sigma$ be a rim hook
contained in the outside border of
$\lambda$. Let
$slitherup(\lambda,\sigma)$ denote
 the rim hook
contained in the outside border
of $\lambda$, whose tail is adjacent
to the head of $\sigma$. Similarly, we
can define
$slitherdown(\lambda,\sigma)$.
Figure \ref{slipe} gives
illustrations of these two operations.

\begin{figure}[h,t]
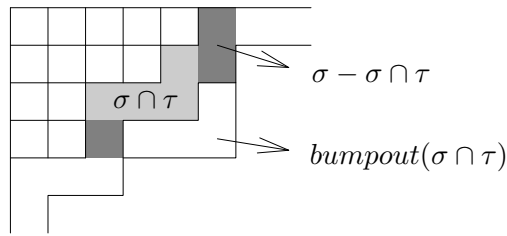

\centertexdraw{ \drawdim mm \linewd 0.1 \setgray 0

\move(10 5)\lvec(25 5)\lvec(25 15)\lvec(20 15)\lvec(20 10)\lvec(10
10)\ifill f:0.8

\move(25 10)\lvec(30 10)\lvec(30 20)\lvec(25 20)\ifill f:0.5

\move(10 0)\lvec(15 0)\lvec(15 5)\lvec(10 5)\ifill f:0.5

\move(15 0)\lvec(30 0)\lvec(30 10)

\textref h:C v:C \htext(18 7.5){$\sigma\cap\tau$}\move(27.5
2.5)\arrowheadtype t:V \avec(37 1)\textref h:C v:C \htext(53
0){$bumpout(\sigma \cap\tau)$}\move(27.5 15)\arrowheadtype t:V
\avec(37 12)\textref h:C v:C \htext(48 11){$\sigma-\sigma\cap\tau$}
 \move(0
0)\lvec(10 0)\move(0 5)\lvec(10 5)\move(0 10)\lvec(20 10) \move(0
15)\lvec(25 15)\move(0 20)\lvec(40 20) \move(0 -10)\lvec(0
20)\move(5 0)\lvec(5 20)\move(10 0)\lvec(10 20)\move(15 10)\lvec(15
20)\move(20 10)\lvec(20 20)\move(25 15)\lvec(25 20) \move(40
15)\lvec(30 15)\lvec(30 10)\lvec(25 10)\lvec(25 5)\lvec(15
5)\lvec(15 -5)\lvec(5 -5)\lvec(5 -10) }
\caption{$\sigma[\tau]$}\label{bump}
\end{figure}

If $\alpha$ is a collection of cells
which are not necessarily contiguous,
then $bumpout(\alpha)$ is the
collection of cells directly bellow and
to the right of a single
cell in $\alpha$. If $\sigma$ and $\tau$ are
two distinct rim hooks outside
$\lambda$ such that $\sigma\cap\tau \neq \emptyset$,
then define $\sigma[\tau]$ as a
 rim hook outside $\lambda\cup \tau$, that is,
\[ \sigma[\tau]=(\sigma-\sigma\cap\tau)\cup
bumpout(\sigma \cap\tau).\]
See Figure \ref{bump} for
an illustration.

Let $U$ be an $m$-rim hook tableau of
shape $\lambda$, and $V$ a
skew $m$-rim hook tableau of shape
$\mu/\omega$. The pair $(U,V)$ is
called an overlapping pair if $\sigma=\lambda/\omega$
is a rim hook
 outside $\omega$, and any content in
 $U$ is smaller than any
content in $V$.

Given an overlapping pair $(U,V)$, Stanton and White
\cite{Stanton} defined
 an operator $A$
acting on the pair $(U,V)$, which generates an
 overlapping pair
$(U^1,V^1)=A(U,V)$. To be precise, the operator $A$
can be described  as follows. Suppose
that $r$ is the minimum content in $V$,
and $\tau$ is the $m$-rim hook
containing $r$. For a rim hook $\alpha$ and
 an integer $s$, let
$\alpha(s)$ denote the rim kook $\alpha$
with  each cell of $\alpha$ filled with $s$.
The skew $m$-rim hook tableau $V^1$ is
obtained from $V$ by removing $\tau(r)$.
 Then
$U^1$ can be constructed
depending on how the shapes $\sigma$ and $\tau$
overlap. Note that $\sigma=\lambda/\omega$,
where $\lambda$ is the
shape of $U$ and $\mu/\omega$ is the
shape of $V$. There are three cases.

Case 1: $\sigma\cap\tau=\emptyset$.
    Set  $U^1=U\cup\tau(r)$.

Case 2:  $\sigma\cap\tau\neq \emptyset$ and
$\sigma\neq\tau$.
Set $U^1=U\cup \tau[\sigma](r)$.

Case 3:   $\sigma=\tau$. In this case,
we construct
a sequence $(\tau=\tau_0, \tau_1, \tau_2,
\ldots)$ of $m$-rim hooks
contained in the outside border
of $\lambda$, where
$\tau_i=slitherup(\lambda,\tau_{i-1})$
 for $i\geq 1$. Assume that  $m_0$ is
the smallest integer such that the
head of $\tau_{m_0}$ is legal with
respect to $\lambda$.
Set $U^1=U\cup\tau_{m_0}(r)$.

It is not difficult to check
that $A(U,V)$ is an overlapping pair,
see \cite{White}.
Let \[A^n(U,V)=A^{n-1}(A(U,V)),\] and write
\[(U^n,V^n)=A^n(U,V).\]
Assume that $n_0$ is the
smallest integer such
that $V^{n_0}$ is
empty. Then, define  \[Combine(U,V)=U^{n_0}.\]

We can now state  the rim
hook insertion algorithm.
Suppose that
$P$ is an $m$-rim hook tableau
 of shape $\lambda$ with contents
 not containing
$r$, and $H$ is an $m$-hook
tableau with content $r$. The rim hook
insertion algorithm gives
 an $m$-rim hook tableau by inserting
$H$ to $P$, denoted $P\leftarrow H$.
To obtain $P\leftarrow H$, we
need to define an overlapping
pair $(U,V)$. Let
$P=P_1\cup P_2$, where $P_1$
(resp., $P_2$) is the $m$-rim hook
tableau contained in $P$ with
contents smaller than (resp., greater
than) $r$. We set $V=P_2$.

To define $U$, we assume that $\tau$ is
the shape of $H$, and $\lambda'$ is the
shape of $P_1$. Let  $ (\tau_0,
\tau_1, \tau_2,
\ldots)$ be a sequence of $m$-rim hooks,
where $\tau_0=\tau$ and
$\tau_i=slitherup(\emptyset,\tau_{i-1})$
 for $i\geq 1$. Moreover, we assume that  $m_0$ is
the smallest integer such that
$\tau_{m_0}$ has no
intersection with $\lambda'$.
We still need another sequence  $(\sigma_0,
\sigma_1,\sigma_2,\ldots)$ of $m$-rim hooks
contained in the outside border
of $\lambda'$, where $\sigma_0=\tau_{m_0}$ and
$\sigma_i=
slitherup(\lambda',\sigma_{i-1})$
for $i\geq 1$. Assume that $n_0$
is the smallest integer such that
$\sigma_{n_0}$ has a legal head
with respect to $\lambda'$. Then $U$
is set to be
$P_1\cup \tau(r)$.

Based on the $m$-rim hook tableaux $U$ and $V$,
the $m$-rim hook tableau  $P\leftarrow H$
can be defined as  $Combine(U,V)$.
It can be shown that the
above rim hook insertion
algorithm is
invertible,  see
\cite{White}.

\begin{figure}[h,t]
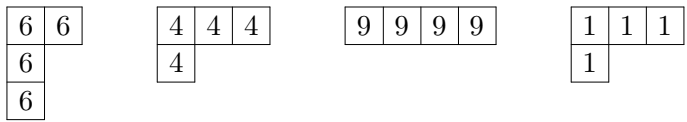

\centertexdraw{ \drawdim mm \linewd 0.1 \setgray 0 \move(0 0)\lvec(5
0)\move(0 5)\lvec(5 5)\move(0 10)\lvec(10 10)\move(0 15)\lvec(10 15)
\move(0 0)\lvec(0 15)\move(5 0)\lvec(5 15)\move(10 10)\lvec(10 15)
\textref h:C v:C \htext(2.5 2.5){$6$}\textref h:C v:C \htext(2.5
7.5){$6$}\textref h:C v:C \htext(2.5 12.5){$6$}\textref h:C v:C
\htext(7.5 12.5){$6$}

\move(20 5)\lvec(25 5)\move(20 10)\lvec(35 10)\move(20 15)\lvec(35
15)\move(20 5)\lvec(20 15)\move(25 5)\lvec(25 15)\move(30
10)\lvec(30 15)\move(35 10)\lvec(35 15)\textref h:C v:C \htext(22.5
7.5){$4$}\textref h:C v:C \htext(22.5 12.5){$4$}\textref h:C v:C
\htext(27.5 12.5){$4$}\textref h:C v:C \htext(32.5 12.5){$4$}

\move(45 10)\lvec(65 10)\move(45 15)\lvec(65 15) \move(45
10)\lvec(45 15)\move(55 10)\lvec(55 15)\move(60 10)\lvec(60 15)
\move(65 10)\lvec(65 15)\move(50 10)\lvec(50 15) \textref h:C v:C
\htext(47.5 12.5){$9$}\textref h:C v:C \htext(52.5 12.5){$9$}
\textref h:C v:C \htext(57.5 12.5){$9$}\textref h:C v:C \htext(62.5
12.5){$9$}

\move(75 5)\lvec(80 5)\move(75 10)\lvec(90 10)\move(75 15)\lvec(90
15)\move(75 5)\lvec(75 15)\move(80 5)\lvec(80 15)\move(85
10)\lvec(85 15)\move(90 10)\lvec(90 15)\textref h:C v:C \htext(77.5
7.5){$1$}\textref h:C v:C \htext(77.5 12.5){$1$}\textref h:C v:C
\htext(82.5 12.5){$1$}\textref h:C v:C \htext(87.5 12.5){$1$}

 } \caption{A $4$-hook permutation}\label{hookpermutation}
\end{figure}

The rim hook insertion algorithm
leads to a one-to-one
correspondence between $m$-hook
 permutations and  pairs of $m$-rim
hook tableaux of the same shape.
 Let $i_1,i_2,\ldots,i_n$ be $n$
positive integers with
$i_1<i_2<\cdots<i_n$. An
$m$-hook permutation
on $\{i_1,i_2,\ldots,i_n\}$
is a permutation of $n$
$m$-hook
tableaux such that the
contents of these $m$-hook
tableaux read off
from left to right form a
permutation on $\{i_1,i_2,
\ldots,i_n\}$. For
example, Figure
\ref{hookpermutation} illustrates
 a $4$-hook permutation
on $\{1,4,6,9\}$.

As remarked in \cite{Stanton}, an $m$-hook
permutation on $[n]$
can be viewed as an element of the
wreath product $C_m\wr S_n$, where  $C_m$ is
the cyclic group of order $m$, or equivalently,
as an $m$-colored permutation
on $[n]$ in the sense that each element
in the permutation is
assigned one of the colors $c_1,c_2,\ldots,c_m$.

Given an $m$-hook permutation
$\mathcal{H}=H_1H_2\cdots H_n$,
one can construct  an $m$-rim hook
tableaux
$P=(\cdots((\emptyset\leftarrow H_1)
\leftarrow H_2)\leftarrow\cdots
)\leftarrow H_n$ by inserting
 $H_1,H_2,\ldots,H_n$ one after another,
see \cite{Stanton} for details.
We call $P$  the insertion
tableau of
$\mathcal{H}$.

\begin{theorem}[\mdseries{Stanton
and White \cite{Stanton}}]\label{ts} There is  a bijection
$\mathrm{Sch}$ between  $m$-hook
permutations on $\{1,2,\ldots,n\}$
and all pairs of $m$-rim hook
tableaux of the same shape with
content $\{1,2,\ldots,n\}$.
\end{theorem}

The bijection $\mathrm{Sch}$ in Theorem \ref{ts}
inherits many
important properties of the ordinary
Schensted correspondence. For our purpose,
 we need the property on the lengths
of the longest increasing and decreasing subsequences
in an $m$-hook
permutation.
An increasing subsequence in an $m$-hook permutation
$\mathcal{H}=H_1H_1\cdots H_n$  is a subsequence
$H_{i_1}H_{i_2}\cdots H_{i_s}$
such that  $H_{i_1},H_{i_2},\ldots,
H_{i_s}$ are of the same shape and
$content(H_{i_1})<content(H_{i_2})<\cdots<content(H_{i_s})$.
A decreasing subsequence in $\mathcal{H}$ can be
 defined analogously.
The following theorem is due to Stanton and White
  \cite{Stanton}, where $\lceil x \rceil$ is
the usual ceiling function meaning the smallest
integer greater than or equal to $x$.

\begin{theorem}\label{L}
Let $\mathcal{H}$ be an $m$-hook permutation, and let  $P$ be the
insertion tableau of $\mathcal{H}$. Suppose that $P$ has $r$ rows
and $c$ columns. Then the length of the longest increasing (resp.,
decreasing) subsequence in $\mathcal{H}$ is $\lceil r/m\rceil$
(resp., $\lceil c/m\rceil$).
\end{theorem}

The following proposition will be used in
the proof of Theorem \ref{main11}.

\begin{proposition}\label{prop}
Let $\mathcal{H}=H_1\cdots H_r\cdots H_n$ be an $m$-hook permutation
with $H_r$ containing the maximum content, and let
$\hat{\mathcal{H}}=H_1\cdots \hat{H_r}\cdots H_n$ be the $m$-hook
permutation obtained from $\mathcal{H}$ by deleting $H_r$. Suppose
that $P_1$ (resp., $P_2$) is the insertion tableau of $\mathcal{H}$
(resp., $\hat{\mathcal{H}}$). Then $P_2$ is the tableau obtained
from $P_1$ by removing the $m$-rim hook filled with the maximum
content.
\end{proposition}

\section{Oscillating rim hook tableaux}

In this section, we present a
bijection between oscillating
$m$-rim hook tableaux
 and $m$-colored matchings. Recall that for an
oscillating $m$-rim hook tableaux  $\lambda$,
$r(\lambda)$ (resp., $c(\lambda)$)
is  the maximum number of
rows (resp., columns) of  shapes
appearing in $\lambda$.

\begin{theorem}\label{main11}
There is a bijection $\phi$
between oscillating $m$-rim hook
tableaux of length $2n$  and $m$-colored matchings
on $[2n]$. Moreover, for any oscillating $m$-rim hook tableau
$\lambda$ we have
\begin{equation}\label{L1}
\lceil r(\lambda)/m
\rceil=\mathrm{ne}(\phi(\lambda))
\end{equation} and
\begin{equation}\label{L2}
\lceil c(\lambda)/m
\rceil=\mathrm{cr}(\phi(\lambda)).\end{equation}
\end{theorem}

\pf We first describe the bijection
 $\phi$ from oscillating
$m$-rim hook tableaux of length $2n$
to $m$-colored
matchings on $[2n]$. Let
\[\lambda=(\lambda^0,\lambda^1,\ldots,\lambda^{2n})\]
be an oscillating $m$-rim hook tableau. We shall recursively define
a sequence
\[(M_0,T_0), (M_1,T_1),\ldots, (M_{2n},T_{2n}),\]where $M_i$ is a
set of $m$-colored arcs, and $T_i$ is an $m$-rim hook tableau of shape
$\lambda^i$. Let $M_0$ be the empty set, and let $T_0$ be the empty
tableau. To obtain $(M_i,T_i)$ for $i\geq 1$, we have the following
two cases.

Case 1: $\lambda^i$ is obtained from $\lambda^{i-1}$
 by adding an
$m$-rim hook outside $\lambda_{i-1}$.
In this case, let
$M_i=M_{i-1}$ and let $T_i$ be the $m$-rim hook
tableau obtained from
$T_{i-1}$ by  filling the $m$-rim hook
$\lambda^i/\lambda^{i-1}$ with the element $i$.

Case 2: $\lambda^i$ is obtained from
$\lambda^{i-1}$ by deleting an
outer $m$-rim hook of
$\lambda_{i-1}$. In this case,
let $T_i$ be
the $m$-rim hook tableau of shape
$\lambda^i$ such that $T_{i-1}$ is
obtained from $T_i$ by inserting
an $m$-hook tableau $H$. Note that
 $H$ is uniquely
determined since the generalized
Schensted
algorithm is invertible. Suppose
that $H$ has shape $(t,1,1,\ldots,1)$
and $content(H)=\{j\}$. It is easy
to see that $j<i$. Let $M_i$
be the set obtained from $M_{i-1}$
by adding an arc $(j, i)$
colored with ${c}_t$.

From the above construction, it is
easy to check  that $M_{2n}$ is an
$m$-colored matching on $[2n]$. We
set $\phi(\lambda)=M_{2n}$.

The reverse map $\psi$, from
$m$-colored matchings to oscillating
$m$-rim hook tableaux, can be described
as follows.
 Let $M$
be an $m$-colored matching on $[2n]$.
We shall construct a sequence
$T(M)=(T_0,T_1,\ldots,T_{2n})$ of
$m$-rim hook tableaux by a
recursive procedure.

Let $T_{2n}$ be the empty tableau.
For $j \leq 2n-1$,
$T_{j}$ can be constructed
based on the following two cases.

\noindent
Case 1: $j+1$ is the right-hand
endpoint of an  arc $(i,j+1)$
colored with ${c}_t$. In this
case, let $T_j$ be the $m$-rim
hook tableau obtained from $T_{j+1}$
by inserting an $m$-hook
tableau of shape $(t,1,\ldots,1)$
filled with $i$.

\noindent
Case 2: $j+1$ is the left-hand
endpoint of an arc $(j+1,k)$.
Note that in this case $j+1$ is
the largest content of $T_{j+1}$.
Let $T_j$ be the $m$-rim
hook tableau obtained from $T_{j+1}$
by deleting the $m$-rim hook
filled with $j+1$.

Set $\lambda=(\lambda^0,\lambda^1,\ldots,
\lambda^{2n})$, where
$\lambda^i$ is the shape of $T_i$.
From the above procedure, it can be seen
 that $\lambda$ is an oscillating
$m$-rim hook tableau.
 Define $\psi(M)=\lambda$.
 It is not hard to check that $\psi$
is the reverse map of $\phi$.

It remains to prove the relations
 \eqref{L1} and \eqref{L2}.
For a $m$-colored matching $M$, let $T(M)=(T_0,T_1,\ldots,T_{2n})$
be the sequence of $m$-rim hook tableaux constructed from $M$. We
now give a recursive procedure
 to generate a sequence of
  $m$-hook permutations
$(\mathcal{H}_1,\mathcal{H}_2,\ldots,
 \mathcal{H}_{2n})$
 from the sequence $T(M)$ of $m$-rim hook tableaux.
Let $\mathcal{H}_{2n}$ be the
empty hook permutation. Suppose that
the $m$-hook permutation $\mathcal{H}_i$
has been constructed. To generate
$\mathcal{H}_{i-1}$,
we consider the following two cases.

\noindent
Case 1: $T_{i-1}$ is obtained
from $T_i$ by inserting an $m$-Hook
tableau $H_j$. In this case, let
$\mathcal{H}_{i-1}=\mathcal{H}_iH_j$.

\noindent
Case 2:  $T_i$ is obtained from
$T_{i-1}$ by adding an $m$-rim hook
filled with   $i$. In this case,
let $\mathcal{H}_{i-1}$ be the hook
permutation obtained from
$\mathcal{H}_i$ by deleting the $m$-hook
tableau with content $i$.

We claim that $T_i$ is the
insertion tableau of $\mathcal{H}_i$.
This can be shown by induction.
Clearly, the claim
 holds for $i=2n$. Suppose
that it is true for $i$,
where $1< i\leq 2n$. We shall
show that the claim holds for $i-1$.
 If $\mathcal{H}_{i-1}=\mathcal{H}_iH_j$, then
the statement is obvious. So it suffices to
consider the case when
$\mathcal{H}_{i-1}$ is the hook
permutation obtained from
$\mathcal{H}_i$ by deleting the $m$-hook
tableau with content $i$. Observe that in
this case, $i$ must be the
largest content appearing in the $m$-hook
tableaux of $\mathcal{H}_i$. In
view of Proposition \ref{prop},
we see that   $T_{i-1}$
is  the insertion tableau
of $\mathcal{H}_{i-1}$.

To finish the proof  of \eqref{L1}, we proceed to show that  $M$ has
a $k$-crossing if and only if there exists  an $m$-hook permutation
$\mathcal{H}_i$   $(1\leq i\leq 2n)$ that contains  a  decreasing
subsequence of length $k$. Suppose that $\mathcal{H}_i=H_{i 1}H_{i
2} \cdots H_{it_i}$ ($1\leq i\leq 2n$) contains a  decreasing
subsequence of
 length $k$. By the construction
 of $\mathcal{H}_i$, there are $t_i$ $m$-colored arcs
 of $M$ whose left-hand endpoints are
 $content(H_{i1}),content(H_{i2}),\ldots, content(H_{it_i})$.
For $1\leq s\leq t_i$, let $w_{is}$ be the right-hand endpoint of
 the arc with left-hand endpoint $content(H_{is})$.
Again, by the construction of the
sequence $(\mathcal{H}_1,\mathcal{H}_2,\ldots,
\mathcal{H}_{2n})$, it can be easily
checked that \begin{equation}\label{L3}
w_{i1}>w_{i2}>\cdots>w_{it_i}.\end{equation}
Denote by  $H_{is_1}H_{is_2}\cdots H_{is_k}$
the decreasing subsequence of $\mathcal{H}_i$.
Clearly, the arcs \[(content(H_{is_1}),w_{is_1}),
\ldots,(content(H_{is_k}),w_{is_k})\]
have the same color.
Moreover, from \eqref{L3} it follows
that they form a
 $k$-crossing of $M$.

Finally, we need to consider  the other direction of the above
statement. Suppose that $M$ has a $k$-crossing consisting the arcs
\[(i_1,j_1), (i_2,j_2),\ldots,(i_k,j_k)\] with
$i_1<i_2<\cdots<i_k<j_1<j_2<\cdots<j_k$. By the construction of
$\phi$, it can be easily checked that $i_1,\ldots,i_k$ are contained
in $content(T_{j_1-1})$. Consider the $m$-hook permutation
$\mathcal{H}_{j_1-1}$. Replacing $i$ with $j_1-1$ in \eqref{L3}, we
obtain that the subsequence of
 $\mathcal{H}_{j_1-1}$
consisting of the $m$-hook tableaux with  contents $i_1,i_2,\ldots, i_k$
is a decreasing subsequence of length $k$. This proves the above claim.
 Thus, relation \eqref{L1} can be deduced from Theorem
\ref{L}.

To prove
\eqref{L2},
it suffices to show that  $M$
has a $k$-nesting
if and only if there exists  an $m$-hook
permutation $\mathcal{H}_i$   $(1\leq i\leq 2n)$
that contains an increasing subsequence of  length $k$.
The proof is omitted  since it is analogous to the case for  $k$-crossings.
 \qed

\begin{figure}[h,t]
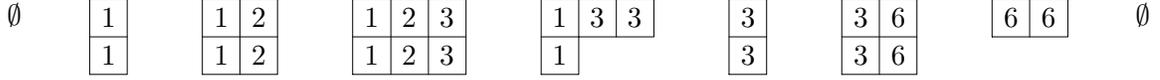

\centertexdraw{ \drawdim mm \linewd 0.1 \setgray 0

\move(10 0)\lvec(15 0)\lvec(15 10)\lvec(10 10)\lvec(10 0)\move(10
5)\lvec(15 5)

\move(25 0)\lvec(35 0)\lvec(35 10)\lvec(25 10)\lvec(25 0)\move(25
5)\lvec(35 5)\move(30 0)\lvec(30 10)

\move(45 0)\lvec(60 0)\lvec(60 10)\lvec(45 10)\lvec(45 0)\move(50
0)\lvec(50 10)\move(55 0)\lvec(55 10)\move(45 5)\lvec(60 5)

\move(70 0)\lvec(75 0)\lvec(75 5)\lvec(85 5)\lvec(85 10)\lvec(70 10)
\lvec(70 0)\move(70 5)\lvec(75 5)\move(75 5)\lvec(75 10)\move(80
5)\lvec(80 10)

\move(95 0)\lvec(100 0)\lvec(100 10)\lvec(95 10)\lvec(95 0)\move(95
5)\lvec(100 5)

\move(110 0)\lvec(120 0)\lvec(120 10)\lvec(110 10)\lvec(110
0)\move(110 5)\lvec(120 5)\move(115 0)\lvec(115 10)

\move(130 5)\lvec(140 5)\lvec(140 10)\lvec(130 10)\lvec(130
5)\move(135 5)\lvec(135 10)

\textref h:C v:C \htext(150 8){$\emptyset$}\textref h:C v:C \htext(0
8){$\emptyset$}\textref h:C v:C \htext(12.5 2.5){$1$}\textref h:C
v:C \htext(12.5 7.5){$1$}\textref h:C v:C \htext(27.5
2.5){$1$}\textref h:C v:C \htext(27.5 7.5){$1$}\textref h:C v:C
\htext(32.5 2.5){$2$}\textref h:C v:C \htext(32.5 7.5){$2$} \textref
h:C v:C \htext(47.5 2.5){$1$}\textref h:C v:C \htext(47.5 7.5){$1$}
\textref h:C v:C \htext(52.5 2.5){$2$}\textref h:C v:C \htext(52.5
7.5){$2$}\textref h:C v:C \htext(57.5 2.5){$3$}\textref h:C v:C
\htext(57.5 7.5){$3$}\textref h:C v:C \htext(132.5 7.5){$6$}\textref
h:C v:C \htext(137.5 7.5){$6$}\textref h:C v:C \htext(112.5
2.5){$3$}\textref h:C v:C \htext(112.5 7.5){$3$}\textref h:C v:C
\htext(117.5 2.5){$6$}\textref h:C v:C \htext(117.5 7.5){$6$}
\textref h:C v:C \htext(97.5 2.5){$3$}\textref h:C v:C \htext(97.5
7.5){$3$}\textref h:C v:C \htext(72.5 2.5){$1$}\textref h:C v:C
\htext(72.5 7.5){$1$}\textref h:C v:C \htext(77.5 7.5){$3$}\textref
h:C v:C \htext(82.5 7.5){$3$}

} \caption{An illustration of the bijection $\phi$}\label{OO}
\end{figure}

Figure \ref{OO} gives an illustration of the bijection $\phi$ for an
oscillating domino tableau.
In this case, the corresponding 2-colored matching is given in
Figure 3.2, where we use solid lines to represent arcs assigned the
color ${c}_1$,
 and dotted lines to
represent arcs with  color ${c}_2$. It can be easily checked
that
\[\lceil r(\lambda)/2 \rceil=\mathrm{ne}(\phi(\lambda))=1\ \ \ \
\text{and}\ \ \ \ \lceil c(\lambda)/2
\rceil=\mathrm{cr}(\phi(\lambda))=2.\]

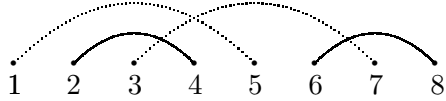
\begin{figure}[h,t]
\setlength{\unitlength}{0.4mm}
\begin{center}
\begin{picture}(140,40)
\put(0,10){\circle*{2}}\put(20,10){\circle*{2}}\put(40,10){\circle*{2}}
\put(60,10){\circle*{2}}\put(80,10){\circle*{2}}\put(100,10){\circle*{2}}
\put(120,10){\circle*{2}}\put(140,10){\circle*{2}}
\qbezier[60](0,10)(40,50)(80,10)\qbezier[60](40,10)(80,50)(120,10)
\qbezier(20,10)(40,30)(60,10)\qbezier(100,10)(120,30)(140,10)

\put(-2,0){1}\put(17.5,0){2}\put(38,0){3}\put(58,0){4}
\put(78,0){5}\put(98,0){6}\put(118,0){7}
\put(-2,0){1}\put(138.5,0){8}
\end{picture}\caption{An example}
\end{center}\label{sdjgol}
\end{figure}

The bijection in Theorem \ref{main11} implies a symmetry property
for the crossing number and the nesting number of $m$-colored
matchings. Let $\lambda=(\lambda^0,\lambda^1, \ldots,\lambda^{2n})$
be an oscillating $m$-rim hook tableau. The conjugate oscillating
$m$-rim hook tableau of $\lambda$, denoted by $\lambda'$, is defined
by $(\mu^0,\mu^1,\ldots, \mu^{2n})$, where $\mu^i$ is the conjugate
of $\lambda^i$. By Theorem \ref{main11}, we get an involution on
$m$-colored matchings satisfying the following symmetry property.

\begin{corollary}\label{coro}
The crossing number $\mathrm{cr}(M)$
and the nesting number
$\mathrm{ne}(M)$ have a symmetric distribution over all  $m$-colored
matchings on $[2n]$, that is, the number of $m$-colored matchings on $[2n]$ with crossing
number $i$ and nesting number $j$ equals the number of $m$-colored matchings on $[2n]$ with
crossing number $j$ and nesting number $i$, for $0\leq i, j \leq n$.
\end{corollary}

We note that the above symmetry property can also be deduced from
the bijection between oscillating
 tableaux and ordinary matchings, see Chen, Deng, Du, Stanley and Yan \cite{Chen}.
More precisely, an $m$-colored matching $M$ can be considered
 as a sequence $(M_1, M_2,\ldots, M_m)$ of disjoint matchings on $[2n]$,
where  $M_i$ is the set of arcs of $M$ colored with $c_i$. It is
obvious that $M$ is $k$-noncrossing if and
 only if each $M_i$ is $k$-noncrossing.
 Nevertheless, Theorem \ref{main11} does not seem to be a direct consequence of the
 correspondence for the case $m=1$.

\section{Oscillating domino tableaux}

In this section, we introduce the structure of packings of Dyck paths
and establish a connection with oscillating domino tableaux with
each shape having at most two columns. By Theorem \ref{main11}, we
see that such oscillating domino tableaux are in one-to-one
correspondence with noncrossing $2$-colored matchings. We show that
oscillating domino tableaux with at most two columns are in
one-to-one correspondence with
 Dyck path packings. This means that there is a correspondence
between Dyck path packings and noncrossing $2$-colored matchings. It
is easy to prove that the number of noncrossing $2$-colored
matchings on $[2n]$ equals $C_nC_{n+1}$, where $C_n$ is the $n$-th
Catalan number.

As will be seen, noncrossing $2$-colored matchings are quite close
to Guy's walks \cite{Guy} in the sense that there is a simple
bijection between these two structures. However, it should be
mentioned that the problem of  counting Dyck path
packings does  not seem to be as easy as counting noncrossing
$2$-colored matchings. Using the identity
\begin{equation}\label{cin}
\sum_{i=0}^n {2n \choose 2i} C_{i} C_{n-j}= C_nC_{n+1},
\end{equation}
 see \cite{Cori}, we obtain the following formula.

\begin{theorem}\label{ncm-th}
For $n\geq 0$, the number of noncrossing $2$-colored matchings on $[2n]$ equals $C_nC_{n+1}$.
\end{theorem}

 Let
$O_{2n,2}$ be  the set of oscillating domino tableaux $\lambda$ of
length $2n$ with each shape appearing in $\lambda$ having at most
two columns. By Theorem \ref{ncm-th}, we see that the set $O_{2n,2}$
is also counted by $C_nC_{n+1}$. We shall give a bijection between
$O_{2n,2}$ and Dyck path packings. A Dyck path of length $2n$ is a
lattice path in the plane from $(0,0)$ to $(2n,0)$ with steps
$(1,1)$ and $(1,-1)$ that never passes below the $x$-axis. A
dispersed Dyck path  is a concatenation of Dyck paths and some
horizontal steps $(1,0)$ on the $x$-axis. In other words, a
dispersed Dyck path  of length $2n$ is a lattice path from $(0,0)$
to $(2n,0)$ with up and  down steps above the $x$-axis and with
horizontal steps on the $x$-axis.   We say that a dispersed Dyck
path $E$ is weakly covered by a Dyck path $D$ if $E$ never goes
above $D$. If $D$ is a Dyck path of length $2n$
 and $E$ is a dispersed Dyck path  of length $2n$
 that is  weakly covered by
$D$, then we say that the pair $(D,E)$ is a
packing of Dyck paths of length $2n$. Denote by $P_{2n}$ the set of Dyck path packings of length $2n$.

\begin{theorem}\label{agf}
There is a bijection between the set $O_{2n,2}$ of oscillating domino tableaux and the set
$P_{2n}$ of Dyck path packings.
\end{theorem}

\pf Let $\mu=(\mu^0,\mu^1,\ldots,\mu^{2n})$ be an oscillating domino
tableau in $O_{2n,2}$, and let
$\lambda=(\lambda^0,\lambda^1,\ldots,\lambda^{2n})$ be the conjugate
oscillating domino tableau of $\mu$, that is, $\lambda^i$ is the conjugate of $\mu^i$
($0\leq i\leq 2n$). Since each shape $\lambda^i$
has at most two rows, we may write $\lambda^i=(u_i,v_i)$, where
$u_i\geq v_i\geq 0$. Clearly, both $u_i+v_i$ and $u_i-v_i$ are
even. For $0\leq i\leq 2n$, let
\[a_i=\frac{u_i+v_i}{2} \quad \mbox{and} \quad
b_i=\frac{u_i-v_i}{2}.\]
Define two lattice paths $D$ and $E$  by setting
\begin{equation}\label{de}
D= ((0,a_0), (1, a_1), \ldots, (2n, a_{2n})) \quad \mbox{and} \quad
E=((0,b_0), (1, b_1), \ldots, (2n, b_{2n})),\end{equation} where a
lattice path is represented by the lattice points.

We proceed to show that  the map $\alpha \colon \mu\longrightarrow (D,E)$ defined by the
above construction is a bijection
from $O_{2n,2}$ to $P_{2n}$. As the first step, we prove that
$(D,E)$ is a packing of Dyck paths of length $2n$. It can be easily checked
that $a_0=a_{2n}=0$ and $(i+1,a_{i+1})-(i,a_i)=(1,1)\ \text{or}\
(1,-1)$.  Hence $D$ is a Dyck path of length $2n$.
Let us consider the possible values of
$(i+1,b_{i+1})-(i,b_i)$. There are two
cases.

\noindent
Case 1:  $\lambda^{i+1}$ is obtained from $\lambda^i$ by adding or
deleting a horizontal domino.  In this case, it is easy to check
that $(i+1,b_{i+1})-(i,b_i)=(1,1)\ \text{or}\ (1,-1)$.

\noindent
Case 2: $\lambda^{i+1}$ is obtained from $\lambda^i$ by adding or
deleting a vertical domino. In this case,  both $\lambda^i$ and
$\lambda^{i+1}$ have two rows with the same number of cells, that
is, $u_i=v_i$ and $u_{i+1}=v_{i+1}$. This implies that
$b_i=b_{i+1}=0$ and $(i+1,b_{i+1})-(i,b_i)=(1,0)$.
So we deduce  that if there is a horizontal step
 on the path $E$, then it lies on the $x$-axis. It is clear that
$b_i\geq 0$. Thus $E$ is a dispersed Dyck path of length
$2n$.

Since $a_i\geq b_i$, we see that $E$ is weakly covered
by $D$. Hence we conclude that $(D,E)$ is a packing of Dyck paths of length $2n$.

The inverse map of $\alpha$ can be described as follows.  Let $(D,E)$ be a Dyck path
packing in $P_{2n}$. Write $D$ and $E$ in the forms as in (\ref{de}).
Set
$\lambda^{i}=(a_i+b_i,a_i-b_i)$, and let $\mu^i$ be the conjugate of
$\lambda^i$. It is easy to check that
$(\mu^0,\mu^1,\ldots,\mu^{2n})$ is an oscillating domino tableau
belonging to $O_{2n,2}$. This completes the proof. \qed

For example, let $\lambda$ be the conjugate of the oscillating domino
tableau  in Figure \ref{OO}. The  corresponding Dyck path
packing is given  in Figure \ref{image}.
\begin{figure}[h,t]
\centertexdraw{ \drawdim mm \linewd 0.15 \setgray 0

\move(0 0)\lvec(10 10)\lvec(20 20)\lvec(30 30)\lvec(40 20)\lvec(50
10)\lvec(60 20)\lvec(70 10)\lvec(80 0)

\move(0 0)\lvec(10 0)\lvec(20 0)\lvec(30 0)\lvec(40 10)\lvec(50
0)\lvec(60 0)\lvec(70 10)

\textref h:C v:C \htext(0 0){\circle*{2.5}}\textref h:C v:C
\htext(10 0){\circle*{2.5}}\textref h:C v:C \htext(10
10){\circle*{2.5}} \textref h:C v:C \htext(20
0){\circle*{2.5}}\textref h:C v:C \htext(20 20){\circle*{2.5}}
\textref h:C v:C \htext(30 0){\circle*{2.5}}\textref h:C v:C
\htext(30 30){\circle*{2.5}}\textref h:C v:C \htext(40
10){\circle*{2.5}}\textref h:C v:C \htext(40 20){\circle*{2.5}}
\textref h:C v:C \htext(50 0){\circle*{2.5}}\textref h:C v:C
\htext(50 10){\circle*{2.5}}\textref h:C v:C \htext(60
0){\circle*{2.5}}\textref h:C v:C \htext(60 20){\circle*{2.5}}
\textref h:C v:C \htext(70 10){\circle*{2.5}}\textref h:C v:C
\htext(80 0){\circle*{2.5}}\move(70.15 10)\lvec(80.15 0)\move(70
9.85)\lvec(79.85 0) } \caption{An example for Theorem
\ref{agf}}\label{image}
\end{figure}

By Theorem \ref{ncm-th} and Theorem \ref{agf}, we  obtain the following relations.

\begin{theorem} For $n\geq 0$,
\begin{equation} |O_{2n,2}|= |P_{2n}|=C_nC_{n+1}.
\end{equation}
\end{theorem}

\begin{figure}[h,t]
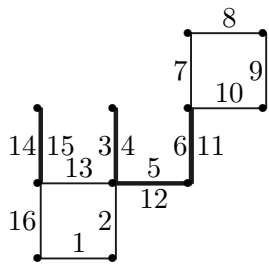

\centertexdraw{ \drawdim mm \linewd 0.22 \setgray 0.1

\move(0 0)\lvec(10 0)\move(0 10)\lvec(20 10)\move(20 20)\lvec(30
20)\move(20 30)\lvec(30 30)

\move(0 0)\lvec(0 10)\move(10 0)\lvec(10 20)\move(20 10)\lvec(20
30)\move(30 20)\lvec(30 30)

\move(0 10)\lvec(0 20)

\textref h:C v:C \htext(0 0){\circle*{2.5}}\textref h:C v:C
\htext(10 0){\circle*{2.5}}\textref h:C v:C \htext(10
10){\circle*{2.5}}\textref h:C v:C \htext(0 10){\circle*{2.5}}
\textref h:C v:C \htext(10 20){\circle*{2.5}}\textref h:C v:C
\htext(20 10){\circle*{2.5}}\textref h:C v:C \htext(20
20){\circle*{2.5}}\textref h:C v:C \htext(20 30){\circle*{2.5}}
\textref h:C v:C \htext(30 20){\circle*{2.5}}\textref h:C v:C
\htext(30 30){\circle*{2.5}}\textref h:C v:C \htext(0
20){\circle*{2.5}}

\move(10 9.8)\lvec(20 9.8)\move(10 10.2)\lvec(20 10.2)

\move(9.8 10)\lvec(9.8 20)\move(10.2 10)\lvec(10.2 20)

\move(19.8 10)\lvec(19.8 20)\move(20.2 10)\lvec(20.2 20)

\move(-0.2 10)\lvec(-0.2 20)\move(0.2 10)\lvec(0.2 20)

\textref h:C v:C \htext(5 2){1}\textref h:C v:C \htext(8.5 5){2}
\textref h:C v:C \htext(8.5 15){3}\textref h:C v:C \htext(11.5
15){4}\textref h:C v:C \htext(15 12){5}\textref h:C v:C \htext(18.5
15){6}\textref h:C v:C \htext(18.5 25){7}\textref h:C v:C \htext(25
32){8}\textref h:C v:C \htext(28.5 25){9}\textref h:C v:C \htext(25
22){10}\textref h:C v:C \htext(22.5 15){11}\textref h:C v:C
\htext(15 8){12}\textref h:C v:C \htext(5 12){13}\textref h:C v:C
\htext(-2.5 15){14}\textref h:C v:C \htext(2.5 15){15}\textref h:C
v:C \htext(-2.5 5){16}

} \caption{A Guy's walk}\label{Guy}
\end{figure}

To conclude this paper, we give a bijection between
noncrossing $2$-colored matchings and Guy's walks.  More precisely, a Guy's walk is
defined to be a
 lattice walk within the  first quadrant starting and ending
at $(0,0)$ and consisting of the following steps
 \[(-1,0),\; (0,-1), \; (0,1), \;
(1,0).\]
Figure \ref{Guy} gives a Guy's walk with steps
\[\rightarrow\ \ \uparrow\ \ \uparrow\ \ \downarrow\ \ \rightarrow\ \ \uparrow\ \  \uparrow
\ \ \rightarrow\ \ \downarrow\ \ \leftarrow\ \ \downarrow\ \ \leftarrow\ \
\leftarrow\ \ \uparrow\ \ \downarrow\ \ \downarrow.
\]

\begin{theorem}\label{cc}
There is a bijection between noncrossing  2-colored matchings on
$[2n]$ and Guy's walks with $2n$ steps.
\end{theorem}

\pf We construct a  bijection $\beta$ from noncrossing 2-colored
matchings on $[2n]$ to Guy's walks with $2n$ steps.  Let $M$ be a noncrossing 2-colored
matching on $[2n]$. For any $i$ in $[2n]$, we define a step
 $s_i$  as follows:
 \[ s_i= \left\{\begin{array}{ll}
       (1,0), \quad & \mbox{if $i$ is the left-hand endoint of an arc with color $c_1$;}\\[6pt]
                  (-1,0), \quad & \mbox{if $i$ is the right-hand endpoint of an arc with color $c_1$;}\\[6pt]
                  (0,1), \quad & \mbox{if $i$ is the left-hand endpoint of an arc with color $c_2$;}\\[6pt]
                   (0,-1), \quad & \mbox{if $i$ is the right-hand endpoint of an arc with  color
${c}_2$.}
\end{array}\right.
\]
Then define $\beta(M)$ to be the walk $s_1s_2\cdots s_{2n}$ starting at the origin.

 We claim  that
$\beta(M)$ is a walk in the first quadrant, starting and ending at the
origin.  For any $i\in [2n]$, consider
the first $i$ steps $s_1,s_2,\ldots,s_i$.  Since
 $M$ is a  noncrossing 2-colored
 matching, it is not hard to check that there are at least as many  $(1,0)$ steps as
 $(-1,0)$ steps  in $\{s_1,s_2,\ldots,s_i\}$. Similarly, there are at least as many
 $(0,1)$  steps as  $(0,-1)$ steps in $\{s_1,s_2,\ldots,s_i\}$.
  Hence $\beta(M)$
 is in the first quadrant.
 It is clear   that $\beta(M)$ terminates  at the
 origin. So the claim holds.  It is not difficult to see that the map   $\beta$ is invertible.
  This completes
the proof.  \qed

\vspace{.2cm} \noindent{\bf Acknowledgments.} This work was
supported by the 973 Project, the PCSIRT Project of the Ministry of
Education, and the National Science Foundation of China.

%------------------------------------------------------------------------

\end{document}